\documentclass[11pt,a4paper]{article}
\usepackage{amsthm,amsmath,amssymb,amsfonts}
\usepackage{graphicx}
\usepackage{url}

\DeclareMathOperator{\K}{\mathrm{K}}

\begin{document}

\newtheorem{proposition}{Proposition}
\newtheorem{lemma}{Lemma}
\newtheorem{theorem}{Theorem}
\newtheorem{corollary}{Corollary}

\def\INT{\mathbb{N}}
\def\ZZ{\mathbb{Z}}
\def\BIT{\mathbb{B}}

\def\eps{\varepsilon}
\def\divides{\mathrel{\vdots}}

\title{Upper bound for the generalized repetition threshold.}

\author{A.\,Rumyantsev\thanks{Moscow State University; \texttt{azrumyan@mail.ru}.
Supported by NAFIT ANR-08-EMER-009-0[1,2] and RFBR 09-01-00709a grants.}}

\date{}

\maketitle

\begin{abstract}
Let $A$ be an $a$-letter alphabet. We consider fractional powers
of $A$-strings: if $x$ is a $n$-letter string, $x^r$ is a prefix
of $xxxx\ldots$ having length $nr$.

Let $l$ be a positive integer. Ilie, Ochem and Shallit defined
$R(a,l)$ as the infimum of reals $r>1$ such that there exist a
sequence of $A$-letters without factors (substrings)
that are fractional powers $x^{r'}$ where $x$ has length at
least $l$ and $r'\ge r$.

We prove that $1+\frac{1}{la}\le R(a,l)\le 1+\frac{c}{la}$ for
some constant~$c$.
\end{abstract}

\section{Introduction}

A fractional power $x^r$ of a string~$x$ is defined as
$x^r=xxx\dots xxy$ where $y$ is a prefix of $x$ and
$|x^r|=r|x|$. (We assume that $r>1$ is a fraction with
denominator $|x|$.)

One may ask whether there exists an infinite sequence of letters
that does not contain fractional powers $x^r$ with large $r$ and
long $x$. More precisely, for a given alphabet size $a$, a given
integer $l$ and a given real $\alpha$ one may ask whether there
exists an infinite sequence of letters that does not contain
fractional powers $x^r$ with $r>\alpha$ and $|x|\ge l$.

For $\alpha=1$ the answer is evidently negative (each string $x$
is a fractional power $x^1$). On the other hand, it is easy to
see that for any $a\ge 2$ and $l\ge 1$ the answer is positive if
$\alpha$ is large enough (there exists a binary sequence that
does not contain factors $x^3$). The threshold value that
separates negative and positive answers is denoted by $R(a,l)$
in~\cite{longpowerfree}; the authors note that $1<R(a,l)\le 2$
and compute exact values of $R(a,l)$ for some pairs~$(a,l)$.
Evidently, $R(a,l)$ decreases when $a$ or $l$ increase.

To get a lower bound for $R(a,l)$, let us apply the pigeonhole
principle to $a+1$ letters at positions $0,l,2l,\ldots,al$. Two
of them should be equal and this creates a fractional power
$x^r$ where $|x|\ge l$ and $r\le 1+1/la$ (this power starts and
ends with a letter that appears twice). Therefore,
    $$
R(a,l)\ge 1+\frac{1}{la}.
    $$

Francesca Fiorenzi, Pascal Ochem and Elise Vaslet
in~\cite{vaslet} gave stronger lower bounds and also some upper
bounds for~$R(a,l)$. In particular, they proved that
    $$
1+\frac{1}{1+\left\lfloor\frac{3l+2}{4}(a-1)\right\rfloor}\le
R(a,l)\le 1+\frac{2\ln l}{l\ln\lambda}+O\left(\frac{1}{l}\right),
    $$
where $\lambda=\frac{(a-1)+\sqrt{(a-1)(a+3)}}{2}$ and a constant
in~$O$ may depend on~$a$ but not on~$l$.

In this paper we use Lov\'asz local lemma to prove a stronger
upper bound for~$R(a,l)$. Our upper bound differs from the lower
bound only by a constant:
    $$
R(a,l)\le 1+\frac{c}{la}
    $$
for some $c$ and for all $a\ge 2$, $l\ge 1$.

\section{Kolmogorov complexity of subsequences}

We present the proof using the notion of \emph{Kolmogorov
complexity} (also called \emph{algorithmic complexity} or
\emph{description complexity}). We refer the reader
to~\cite{livitan} or~\cite{shen-lecture-notes} for the
definition and basic properties of Kolmogorov complexity.

For an infinite sequence~$\omega$ and finite set $X\subset
\mathbb{N}$ let $\omega(X)$ be a string of length $\#X$ formed
by $\omega_i$ with $i\in X$ (in the same order as in $\omega$).

We use the following result from~\cite{csr-2007} that guarantees
the existence of a sequence $\omega$ such that strings
$\omega(X)$ have high Kolmogorov complexity for all simple~$X$:

\begin{theorem}\label{subsequence1}
Let $\alpha$ be a positive real number less than~$1$. There
exists a binary sequence $\omega$ and an integer $N$ such that
for any finite set $X$ of cardinality at least $N$ the
inequality
    $$
\K(X,\omega(X)|t)\ge\alpha\#X
    $$
holds for some $t\in A$.
\end{theorem}

Here $\K(X,\omega(X)|t)$ is conditional Kolmogorov complexity
of a pair $(X,\omega(X))$ relative to~$t$.

We need a slightly more general version of this result (for any
alphabet size):

\begin{theorem}\label{subsequencea}
Let $a\ge 2$ be an integer. Let $\alpha$ be a positive real less
than~$1$. There exists a sequence $\omega$ in $a$-letters
alphabet and an integer $N$ such that for any finite set $X$ of
cardinality at least $N$ the inequality
    $$
\K(X,\omega(X)|t)\ge\alpha\#X\log a
    $$
holds for some $t\in X$.
\end{theorem}

\begin{proof}
Theorem~\ref{subsequencea} can be proven using exactly
the same argument as in~\cite{csr-2007} (Lovasz local lemma
technique). It can also be formally derived from
Theorem~\ref{subsequence1} as follows: we encode $a$ letters of
the alphabet by bit blocks of some length $t$ (large enough).
This encoding is not bijective (several blocks encode the same
letter) but is chosen in such a way that all letters have almost
the same number of encodings (about $2^t/a$). Then we take a
sequence from Theorem~\ref{subsequence1}, split it into $t$-bit
blocks and replace these blocks by corresponding letters. If some
subsequence formed by the letters is simple, then the
corresponding bit subsequence is simple, too. (Technically we
should change $\alpha$ slightly to compensate for ``boundary
effects''.)
\end{proof}

\section{Weak upper bound}

To illustrate the technique, we first prove a simple
generalization of a result obtained by
Berk~\cite{repetitionfree} and provide an upper bound for
$R(a,l)$ that is weaker that our final bound:

\begin{theorem}
For every $a\ge 2$ and every real number $b\in(1,a)$ there exists a
number~$N$ and a sequence~$\omega$ in $a$-letters alphabet such
that for every $n\ge N$ the distance between any two different
occurrences of the same substring of length~$n$ in $\omega$ is at
least~$b^n$.
\end{theorem}

\begin{proof}
Construct a sequence~$\omega$ using Theorem~\ref{subsequencea}
with $\alpha$ close enough to~$1$.

Let $I$ and $J$ ($|I|=|J|=n$) be different intervals where the
same substring of length $n$ occurs in~$\omega$. Let $X=I\cup
J$. Then $n<\#X\le 2n$ (intervals $I$ and $J$ are not
necessarily disjoint) and the first $n$ letters of $\omega(X)$
are equal to the last $n$ letters of $\omega(X)$. It is easy
to see that the string $\omega(X)$ is determined
by its first $\#X-n$ letters, $n$ and $\#X$, so
$\K(\omega(X))\le(\#X-n)\log a+O(\log n)$.

Assume $t\in X$. Then $X$ is determined by $t$, the number $n$,
the distance between $I$ and $J$ and the ordinal number of $t$
in $X$. So if the distance between $I$ and $J$ is less than
$b^n$ then $\K(\omega(X),X|t)\le(|X|-n)\log a+ n\log b+O(\log
n)\le\alpha|X|\log n$ for large enough~$n$ and $\alpha$ that is
close enough to~$1$ (because $\log b<\log a$). This contradicts
the inequality of Theorem~\ref{subsequencea}. Therefore
sequence~$\omega$ does not contain a pair of different
occurrences of the same substring of sufficiently large length~$n$
with distance between them less than $b^n$.
\end{proof}

In particular, for every integer $a\ge 2$, every real number
$b\in(1,a)$ and for large enough $l$ the following inequality
holds:
    $$
R(a,l)<1+\frac{\log_b l}{l}.
    $$

\section{The final upper bound}

In the weak upper bound we used the same sequence for all values
of~$l$. And now we need different sequences for different values
of~$l$ but we want the constant~$c$ to be the same. To achieve
this goal we use the following ``$l$-uniform'' version of
Theorem~\ref{subsequence1}.

\begin{theorem}\label{subsequencel}
Let $\alpha$ be a positive real number less than~$1$. There
exists an integer $N$ such that for every integer~$l$ there
exists a binary sequence $\omega$ that has the following
property: for every finite set $X$ of cardinality at least $N$
the inequality
    $$
\K(X,\omega(X)|t,l)\ge\alpha\#X
    $$
holds for some $t\in A$.
\end{theorem}

Note that $\omega$ may depend on~$l$ while $N$ is the same for
all values of~$l$. (If we allowed $N$ to be dependent on $l$,
this would be a standard relativization of
Theorem~\ref{subsequence1}.)

\begin{proof}
Theorem~\ref{subsequencel} can be proven in the same way as
Theorem~\ref{subsequence1}. And it can also be formally derived
from it: if a sequence~$\tau$ and a number~$N$ satisfy the
requirements of Theorem~\ref{subsequence1} and $z:\INT^2\to\INT$
is a computable bijection, then the sequence $i\mapsto
\omega_i=\tau_{z(i,l)}$ and the same number~$N$ satisfy the
requirements of Theorem~\ref{subsequencel} for the integer~$l$.
(The bijection adds $O(1)$-term, but this can be compensated by
a small change in $\alpha$: the statement is true for every
$\alpha<1$.)
\end{proof}

Now we can start proving the upper bound.

\begin{theorem}
There exists a constant~$c$ such that for any $a\ge 2$ and $l\ge 1$ the
following inequality holds:
    $$
1+\frac{1}{al}\le R(a,l)\le 1+\frac{c}{al}.
    $$
\end{theorem}

\begin{proof}
The lower bound is easy (as shown in the introduction). Let us
prove the upper bound. Let as assume first that $a=2$ (the
general case can be reduced to this special one).

Consider a sequence~$\omega$ satisfying the requirements of
Theorem~\ref{subsequencel} for some $\alpha>\frac{1}{2}$. Then
the required sequence with long fractional powers will be
constructed as
    $$
\tau_i=\omega_{f(i)}
    $$
for some mapping $f:\INT\to\INT$.

At first let us define $f$ at the first $l$ integers (the value
of integer constant~$m$ will be chosen later):
\begin{itemize}
\item[1.] $f(i)=i\bmod m$ for $i<l$ and $(i\bmod m)\neq m-1$ (we say that
these indexes have rank~$1$).
\item[2.] $f(mi+m-1)=(m-1)+(i\bmod m)$ for $mi+m-1<l$ and $(i\bmod m)\neq m-1$
(we say that these indexes have rank~$2$).
\item[3.] $f(m^2i+m^2-1)=2(m-1)+(i\bmod m)$ for $m^2i+m^2-1<l$ and $(i\bmod m)\neq m-1$
(we say that these indexes have rank~$3$).
\item[] (And so on until $f$ is defined at all first $l$ integers.)
\end{itemize}

Then we define $f$ on other blocks of $l$ integers in the same
way but using fresh bits each time. So if
$f(\{0,1,\dots,l-1\})=\{0,1,\dots,L-1\}$ then $f(i+jl)=f(i)+jL$.

Suppose the sequence $\tau_i=\omega_{f(i)}$ contains some
fractional power $xyx$ with $|xy|\ge l$ and the exponent
$\displaystyle \frac{|xyx|}{|xy|}\ge 1+\frac{c}{2l}$. Without
loss of generality we can assume that the exponent
$1+\frac{c}{2l}$ is not greater than~$2$ (otherwise the
statement of the theorem follows from the existence of a binary
sequence, called Thue-Morse sequence, that does not contain any
fractional power with exponent greater than~$2$,
see~\cite{thue-1},~\cite{thue-2}). Also we can assume that
$c>2m$ (increasing $c$, we make our task easier). So
$l\ge\frac{c}{2}>m$ and $|x|\ge\frac{c}{2l}|xy|>m$.

First we consider the case when both occurrences of~$x$ in $xyx$
lie entirely in some blocks of size~$l$ (in two different
blocks, because $|xy|\ge l$). Denote by $n$ the number of $l$-sized
blocks between these two occurrences of $x$ and denote by $k$ the integer
number that satisfies the inequality $m^{k-1}\le|x|<m^k$.
Then $m^k>\frac{c}{2}n$ and $k\ge 2$ (because $|x|\ge\frac{c}{2l}|xy|>m$).

Let us denote by $I$ and $J$ the sets of values of $f$ for the
first and second occurrences of $x$ (respectively) whose rank is
not greater than~$k$ (obviously there is at most $1$ index in
each of these occurrences of~$x$ whose rank is greater
than~$k$). The sets $I$ and $J$ are disjoint because these
occurrences of~$x$ lies in the different $l$-sized blocks.
Assume $Z=I\cup J$, then for some $t\in Z$ we have
$\K(Z,\omega(Z)|t,l)\ge\alpha\#Z$ by the statement of
Theorem~\ref{subsequencel} (we need here that $m>N+1$ since
$\#Z$ should be greater than~$N$).

Obviously,
    $$
\frac{1}{2}\#Z=\#I+O(1)=\#J+O(1)=(k-1)(m-1)+\frac{|x|}{m^{k-1}}+O(1).
    $$
The set $Z$ is determined by $t$, $l$, $m$, $n$, $k$, $|x|$ and
the start/end positions for the two occurrences of the word~$x$
modulo $m^k$ (and one bit saying whether $t$ belongs to the
first occurrence of $x$ or to the second one). So $\K(Z\mid
t,l)\le\log n+O(\log(m^k))=O(k\log m)$ (since
$m^k>\frac{c}{2}n$). We can also calculate $\omega(Z)$ if
$\omega(I)$ is given (we need at most one extra bit for
calculating the entire string~$x$). Therefore
    $$
O(k\log m)+\frac{1}{2}\#Z\ge\alpha\#Z,
    $$
but $\alpha>\frac{1}{2}$ and $\#Z\ge 2(k-1)(m-1)+O(1)\ge
k(m-1)+O(1)$. So $k(m-1)<O(k\log m)$ that is a contradiction if
$m$ is large enough. (Recall that the choice of $m$ was
postponed.)

Consider now the general case for the position of the two occurrences
of~$x$. If length of $x$ is not large, i.e. $|x|\le l$, we can reduce
this case to the previous one by splitting $x$ into parts and choosing the
largest part (we must multiply the constant~$c$ by $3$).
Now let $x$ be longer than the block size ($|x|>l$). We can assume that
there is no $l$-sized block that intersects both occurrences of $x$
(in the other case we also split the word~$x$ in parts).

Let us denote by $I$ and $J$ the sets of values of $f$ in the
first and second occurrences of $x$ respectively. The sets
$I$ and $J$ are disjoint. Assume $Z=I\cup J$. Then for some
$t\in Z$ we have $\K(Z,\omega(Z)|t,l)\ge\alpha\#Z$.

The set $Z$ is determined by $t$, $l$, $m$ and the relative
start/end positions of the two occurrence of the word~$x$ with
respect to the one of the preimages of $t$ (for example, the
first one). So $\K(Z\mid t,l)\le\log|xy|+O(\log l)=O(\log|x|)$
(since $|x|\ge l$ and $|x|\ge\frac{c}{2l}|xy|$). To compute
$\omega(Z)$, it is enough to know at most a half of it
($\omega(I)$ or $\omega(J)$, whichever is smaller). Therefore
    $$
O(\log|x|)+\frac{1}{2}\#Z\ge\alpha\#Z,
    $$
but $\alpha>\frac{1}{2}$ and
$\#Z=\Omega\left(\frac{|x|}{l}(m-1)\log_m l\right)=
\Omega\left((\log|x|)\frac{m-1}{\log m}\right)$
(here we use that $|x|>l>m$ and $\frac{|x|}{\log|x|}\ge\frac{l}{\log l}$).
That is a contradiction if $m$ is large enough.

This finishes the proof for $a=2$.

Assume now that $a\ge 6$ and $a$ is even. Let $\omega$ be the
sequence constructed for binary alphabet and
$l'=\frac{a-2}{2}l$. To get the required sequence~$\nu$ we will
color the terms of~$\omega$ into $\frac{a}{2}$ colors: the
$i$-th block of size~$l$ gets color $i\bmod\frac{a}{2}$. Then
the size of the alphabet of sequence~$\nu$ (whose terms are now
$\langle$bit, color$\rangle$ pairs) equals to $a$ and $\nu$ does
not contain fractional powers $z^p$ with $|z|\ge\frac{a-2}{2}l$
and $p\ge 1+\frac{c}{(a-2)l}$. And obviously $\nu$ does not
contain any fractional powers $z^p$ with
$l\le|z|\le\frac{a-2}{2}l$ (because it does not contain pairs of
equal letters at these distances).

Therefore $R(a,l)\le 1+\frac{c}{(a-2)l}$ if $a\ge 6$ and $a$ is
even, and $R(2,l)\le 1+\frac{c}{2l}$.

To prove the theorem for arbitrary $a$ it remains to note that
that $R(a,l)$ is decreasing in $a$, so $R(a,l)\le
1+\frac{3c}{al}$ for every $a\ge 2$, $l\ge 1$.
\end{proof}

\section{Acknowledgements}

The author is grateful to Gregory Kucherov who explained this
problem to the author and suggested to apply the Kolmogorov
complexity technique to it, and to Anna Frid who encouraged the
author to write down the proofs.


\begin{thebibliography}{11}

\bibitem{livitan}
Li~M., Vitanyi~P, \textit{An Introduction to Kolmogorov Complexity
and Its Applications}, 2nd ed. N.Y.: Springer, 1997.

\bibitem{thue-1}
Axel Thue, \textit{\"Uber unendliche Zeichenreihen}, Norske Vid. Skrifter
I Mat.-Nat. Kl., Christiania 7 (1906) 1--22.

\bibitem{thue-2}
Axel Thue, \textit{\"Uber die gegenseitige Lage gleicher Teile gewisser
Zeichenreihen}, Norske Vid. Skrifter I Mat.-Nat. Kl.,
Christiania 1 (1912) 1--67.

\bibitem{cubefreepartial}
Florin Manea, Robert Merca, \textit{Freeness of partial words},
Theoretical Computer Science, 389, Issue 1-2 (December 2007), pp.~265--277.

\bibitem{squarefreepartial}
Vesa Halava, Tero Harju and Tomi K\"arki, \textit{Square-free partial words},
Information Processing Letters, Volume 108, Issue 5 (15 November 2008), pp.~290--292.

\bibitem{repetitionfree}
        %
J.~Berk, An application of Lov\'asz local lemma: there exists an
infinite 01-sequence containing no near identical intervals.
        %
In: A.~Hajnal, L.~Lov\'asz, and V.~T.~S\'os, editors,
\textit{Finite and Finite Sets}, Vol.~37 of Colloq. Math.
Soc. J\'anos Bolyai, 1981, pp.~103--107.

\bibitem{longpowerfree}
        %
Lucian Ilie, Pascal Ochem, Jeffrey Shallit, A Generalization of
Repetition Threshold,
        %
\emph{Mathematical foundations of computer science},
\textbf{345}, Issue 2-3 (November 2005), pp.~359--369.

\bibitem{vaslet}
Francesca Fiorenzi, Pascal Ochem, Elise Vaslet,
Bounds for the generalized repetition threshold.
        %
In: \emph{Theoretical Computer Science, 2010}, submitted, available from
\url{http://www.lri.fr/~fiorenzi/Publications/FiorenziOchemVaslet2010.pdf}.

\bibitem{csr-2007}
Andrey Yu.~Rumyantsev, Kolmogorov Complexity,
Lov\'asz Local Lemma and Critical Exponents,
        %
In: \emph{Computer Science in Russia, 2007}. Lecture Notes in
Computer Science, Volume 4649, Springer, 2007, pp.~349--355.

\bibitem{shen-lecture-notes}
        %
Alexander Shen, Algorithmic Information Theory and Kolmogorov
Complexity, December 2000, lecture notes. Published as
Technical Report 2000-034, Uppsala University,
\url{http://www.it.uu.se/research/publications/reports/2000-034}.

\end{thebibliography}
\end{document}